\documentclass{amsart}
\usepackage{fullpage}
\usepackage{amsmath}
\usepackage{amssymb}
\begin{document}
\title[LINEAR STABILITY OF ... -- JDSGT Vol. 4, Number 1 (2006)]{LINEAR STABILITY OF TRIANGULAR EQUILIBRIUM POINTS
IN THE GENERALIZED PHOTOGRAVITATIONAL RESTRICTED THREE BODY
PROBLEM WITH POYNTING-ROBERTSON DRAG \footnote {\it
\\Journal of Dynamical
Systems $\&$ Geometric Theories \\ Vol. 4, Number 1 (2006) 0-0.\\
\copyright Taru Publications}}
\author[B. ISHWAR and
B.S. KUSHVAH]{{\bf B. ISHWAR$^a$ and
B.S. KUSHVAH$^b$}\\
$a$ P.I  DST Project,$b$ JRF  DST Project\\
    University department of mathematics\\
    B.R.A. Bihar University Muzaffarpur-842001\\
    Bihar (India)\\
   Email:ishwar\_bhola@hotmail.com, bskush@hotmail.com }
\maketitle
\begin{center}
\tiny (Received 27 October 2005)
\end{center}
\begin{abstract}In this paper we have examined the linear stability of triangular
equilibrium points in the generalised photogravitational restricted three body problem with Poynting-Robertson drag. We
have found the position of triangular equilibrium points of our problem. The problem is generalised in the sense that smaller
primary is supposed to be an oblate spheroid.The bigger primary is considered  as radiating. The equations of motion are affected by
radiation pressure force, oblateness and P-R drag. All classical results involving photogravitational and oblateness in restricted
three body problem may be verified from this result. With the help of characteristic equation, we discussed the stability. Finally we
conclude that triangular equilibrium points are unstable.\end{abstract}
 \noindent
  


\section{Introduction}

We study the motion of three finite bodies in the three body problem. The
problem is restricted in the sense that one of the three masses is taken to be so
small that the gravitational effect on the other masses by third mass is
negligible. The smaller body is known as infinitesimal mass and remaining two
massive bodies as finite masses. The classical restricted three body problem is
generalised to include the force of radiation pressure, the Poynting-Robertson
effect and oblateness effect. J.H. Poynting(1903) considered the effect of the
absorption and subsequent re-emission of sunlight by small isolated particles in
the solar system. His work was later modified by H.P. Robertson(1937) who used a
precise relativistic treatments of the first order in the ratio of the velocity
of the particle to that of light.The motion of an infinitesimal body near one of
the Lagrangian points is said to be stable if given a very small displacement
and small velocity, the body oscillates for a considerable time around the point.
If it departs from the point as the time increases, the motion is said to be
unstable. Mathematically, we express that the stability of motion means bounded
displacement and velocity with bounded functions of time in the neighborhood of
the equilibrium points.\\

If a satellite flies high enough above the Earth and is large enough in size, but at the same time has a sufficiently small mass, then the pressure of the Sun's light rays (known as radiation pressure ) has a
very strong effect on its motion. The distance of the satellite to the Sun
practically does not change, and so the magnitude of the radiation pressure
force is practically constant. The Poynting-Robertson (P-R) effect is the most
important nongravitational and nondisipative  interplanetary dust particles. The
P-R effect may be significant in dealing with other problems in astrophysics.\\

The effect of radiation pressure and P-R drag in the
restricted three body problem has been studied by Colombo {\it et al.}(1966),Chernikov Yu.A.(1970) and Schuerman(1980)  discussed the position as well as the stability of the Lagrangian equilibrium points when radiation pressure, P-R drag force are included. Murray C.D.(1994) systematically discussed the dynamical effect of general drag in the planar circular restricted three body problem. Liou J.C. {\it et al.}(1995) examined the effect of radiation pressure, P-R drag and solar wind drag in the restricted three body problem.The classical three body problem has five Lagrangian points. Their location and stability properties are well known. The three collinear points $L_1,L_2,L_3$ lies on the line joining the two primaries, are linearly unstable. The two Lagrangian points $L_4,L_5$ lie in the orbit of the test particle which
are symmetrical about the line joining the primaries, are linearly stable with the condition that the mass parameter $\mu$ to be less than .0385. Hence, we are interested to discuss the linear stability of triangular equilibrium points in generalised photogravitational restricted three body problem with Poynting-Robertson drag. After Newton solved the problem of the orbit of single planet around the Sun, the natural next challenge was to find the solution of two planets orbiting the Sun(i.e. three body problem ). Many of the best minds in Mathematics and Physics worked on this problem in the last century. Three-Body problem is a continuous source of study, since the discovery of its non-integrability due to Poincare(1892). Regular and chaotic motion have been widely investigated with any kind of tools, from analytical result to numerical explorations. The restricted three-body problem is a modified three-body problem when one of the three masses is set to zero(infinitesimal mass ).\\

The restricted three body model, Earth-Moon and Spacecraft is insoluble, although specific solutions exist, like the ones in which the spacecraft is positioned at one of the  Lagrangian points. There are many periodic orbits in the restricted three body problem. One of  the most famous, discovered by Lagrange, is formed by an equilateral triangle. In \lq\lq Earth-Moon-Space Station  model'', if the Moon is not too massive,  this orbit is thought to be stable. If the space station is pushed a bit to one  side (in position or velocity ), it's supposed to make small oscillations around  this orbit. The two kinds of triangular points are called $L_4$ and $L_5$ points. \\

 There is a group of enthusiasts who want to setup a colony at $L_5$ point of   the Earth-Moon system. As already noted, because $L_4$ and $L_5$ are the stable points of equilibrium, they have been proposed for sites of large self-contained  \lq\lq Space colonies'', an idea developed
and advocated by the late Geraid O' Neill (1974). The three body problem have an interesting application for artificial satellites and future space colonization. Triangular points of the Sun-Jupiter or Sun-Earth system would be convenient sites to locate future space colonies. Such colonies would presumably mine asteroids for at least some of their raw materials. Present day technology permits the construction of spacecraft whose value of $\beta$(ratio of radiation force and gravitation force) can be made or varied in the range $0<\beta<1$. (A spacecraft with $\beta>1$ has been considered  by NASA as a  \lq\lq Sailing'' vehicle for  chasing Halley's comet [Friedman 1976]). A Space Colony with a solar facing adjustable \lq\lq sail" could thus \lq\lq park"(i.e. remain stationary with respect to the Sun and a planet without expending large amounts of energy) not only at the classical positions but at any heleocentricdistance out to that of the planet.\\

Thus we have included the purpose of the work and its significance in the field of Celestial Mechanics and Dynamical Astronomy. Application of results to realistic actual problem is obvious. $L_4$ and $L_5$ points became unstable due to P-R drag which is very remarkable and important, where as they are stable in classical result. 

\section{Location of Triangular Equilibrium Points}
Equations of motion are,
\newcommand{\abc}{{\Bigl(1+\frac{5}{2}A_2\Bigr) }}
\newcommand{\aac}{{(1-\frac{A_2}{2}) }}
\newcommand{\amc}{{(1-\mu) }}
\newcommand{\adc}{{\frac{\delta ^2}{2} }}
\newcommand{\zab}{{\Bigl(1+\frac{5A_2}{2r^2_2}\Bigr)}}
\newcommand{\zw}{{\frac{W_1}{r^2_1}}}
\newcommand{\zwf}{{\frac{W_1}{r^4_1}}}
\newcommand{\zx}{{(x+\mu)}}
\newcommand{\zox}{{(x+\mu-1)}}
\newcommand{\zd}{{\displaystyle}}
\newcommand{\zabs}{{\Bigl(1+\frac{5A_2}{2{r^2_2}_*}\Bigr)}}
\newcommand{\zws}{{\frac{W_1}{{r^2_1}_*}}}
\newcommand{\zwfs}{{\frac{W_1}{{r^4_1}_*}}}
\newcommand{\zxs}{{(x_*+\mu)}}
\newcommand{\zoxs}{{(x_*+\mu-1)}}
\begin{align}
\ddot{x}-2n\dot{y}&=U_x ,\quad\text{where},\quad U_x=\frac{\partial{U_1}}{\partial{x}}-\frac{W_{1}N_1}{r^2_1}\\
\ddot{y}+2n\dot{x}&=U_y,\hspace{.85in}U_y=\frac{\partial{U_1}}{\partial{y}}-\frac{W_{1}N_2}{r^2_1}\\
\ddot{z}&=Uz,\hspace{.85in} U_z=\frac{\partial{U_1}}{\partial{z}}-\frac{W_{1}N_3}{r^2_1}\\
U_1&=\zd{\frac{n^2(x^2+y^2)}{2}}+\frac{\amc{q_1}}{r_1}+\frac{\mu}{r_2}+\frac{\mu{A_2}}{2r^3_2}
\end{align}
\begin{gather*}
r^2_1=\zx^2+y^2,\quad  r^2_2=\zox^2+y^2,\quad n^2=1+\frac{3}{2}A_2,\\
N_1=\frac{\zx{N}}{r^2_1}+\dot{x}-ny,\quad N_2=\frac{yN}{r^2_1}+\dot{y}+n\zx,\quad  N_3=\frac{zN}{r^2_1}+\dot{z},\\
 N=\zx\dot{x}+y\dot{y}+z\dot{z}
\end{gather*}

$A_2=\frac{r^2_e-r^2_p}{5r^2}$ be the oblateness coefficient, $r_e$ and$r_p$ be
the equatorial and polar radii respectively $r$ be the distance between
primaries,  $q=\bigl(1-\frac{F_p}{F_g}\bigr)$ be the mass reduction factor
expressed in terms of the particle's radius $a$, density $\rho$ and radiation
pressure efficiency factor $\chi$ (in the C.G.S.system) i.e.,
$q=1-\zd{\frac{5.6\times{10^{-5}}\chi}{a\rho}}$. Assumption $q=constant$ is
equivalent to neglecting fluctuation in the beam of solar radiation and the
effect of solar radiation, the effect of the planet's shadow, obviously
$q\leq1$. Triangular equilibrium points are given by $U_x=0,U_y=0,z=0,y\neq{0}$, then we have
\begin{eqnarray}\label{eq:xstar} x_*&=&x_0\Biggl\{1-\zd{\frac{nW_1\bigl[\amc\abc+\mu\aac\adc\bigr]}{3\mu\amc{y_0 x_0}}}-\adc\frac{A_2}{x_0}\Biggr\} \\\label{eq:ystar}
y_*&=&y_0{\Biggl\{1-\zd{\frac{nW_1\delta^2\bigl[2\mu-1-\mu(1-\frac{3A_2}{2})\adc+7\amc\frac{A_2}{2}\bigr]}{3\mu\amc{y^3_0}}}-\zd{\frac{\delta^2\bigl(1-\adc)A_2}{y^2_0}}\Biggr\}^{1/2} }
\end{eqnarray}
where $x_0=\adc-\mu$,$y_0=\pm\delta\bigl(1-\frac{\delta^2}{4}\bigr)^{1/2}$ and $\delta=q^{1/3}_1$

where subscript * indicate the value at the equilibrium points. Eq.(~\ref{eq:xstar}) and (~\ref{eq:ystar}) are valid for $W_1\ll1$, $A_2\ll1$.

\section{Stability of Triangular Equilibrium Points}

\quad We suppose $\alpha=A_1e^{\lambda{t}}$,$\beta=B_1e^{\lambda{t}}$,$\gamma=C_1e^{\lambda{t}}$ be the small displacements
from Lagrangian points$(x_*,y_*,z_*)$ and $A_1,B_1,C_1$, $\lambda$ are
parameters, $x=x_*+\alpha$,\quad $y=y_*+\beta$ \quad  \text{and} \quad  $z=z_*+\gamma$.
 Therefore the  equations of perturbed motion corresponding to the system of Eq.(1),(2),(3) may be written as,
\begin{align}
\ddot{\alpha}-2n\dot{\beta} &= {U^0_x}+{\alpha}{U^0_{xx}}+{\beta}{U^0_{xy}}+
{\gamma}{U^0_{xz}}+\dot{\alpha}{U^0_{x\dot{x}}}+{\dot{\beta}}{U^0_{x\dot{y}}}+\dot{\gamma}{U^0_{x\dot{z}}} \\
\ddot{\beta}+2n\dot{\alpha}&= {U^0_y}+{\alpha}{U^0_{yx}}+{\beta}{U^0_{yy}}+{\gamma}{U^0_{yz}}+\dot{\alpha}{U^0_{y\dot{x}}}+\dot{\beta}{U^0_{y\dot{y}}}+\dot{\gamma}{U^0_{y\dot{z}}}\\
\ddot{\gamma}&={U^0_z}+{\alpha}{U^0_{zx}}+{\beta}{U^0_{zy}}+{\gamma}{U^0_{zz}}+\dot{\alpha}{U^0_{z\dot{x}}}+\dot{\beta}{U^0_{z\dot{y}}}+\dot{\gamma}{U^0_{z\dot{z}}}
\end{align}
 \[
\text{At} \quad (x_*,y_*,z_*) \quad  U^0_x=U^0_y=U^0_z=U^0_{xz}=U^0_{yz}=U^0_{zx}=U^0_{zy}=U^0_{x\dot{z}}=U^0_{y\dot{z}}=U^0_{z\dot{x}}=U^0_{y\dot{z}}=0
\]
Above system of equations can be written as,

\begin{align}\label{eq:lambda_x}
(\lambda^2-\lambda{U^0_{x\dot{x}}}-{U^0_{xx}})
+[-(2n+{U^0_{x\dot{y}}})\lambda-{U^0_{xy}}]B_1&=0\\\label{eq:lambda_y}
[(2n-{U^0_{y\dot{x}}})\lambda-{U^0_{yx}}]A_1
+(\lambda^2-\lambda{U^0_{y\dot{y}}}-{U^0_{yy}})B_1&=0\\
(\lambda^2-\lambda{U^0_{z\dot{z}}}-U^0_{zz})C_1&=0 \label{eq:z_eq}
\end{align}
Eq.(~\ref{eq:z_eq}) gives us,
$\lambda=\zd{\frac{U^0_{z\dot{z}}\pm\sqrt{(U^0_{z\dot{z}})^2+4U^0_{zz}}}{2}}$. From Eq.(3)
\[
U^0_{zz}=-f,\quad U^0_{z\dot{z}}=-\zw,\quad \text{where,}\quad f=\frac{\amc{q_1}}{r^3_1}+\zd{\frac{\mu}{r^3_2}\biggl(1+\frac{3}{2}\frac{A_2}{r^2_2}\biggr)}
\]
\begin{equation}
\lambda=\zd{\frac{1}{2}\Biggl\{-\zw \pm\sqrt{\frac{W^2_1}{r^4_1}-4f}\Biggr\}}\quad\text{as}\quad W_1\geq 0 \quad \text{and}\quad f>\frac{W^2_1}{4r^4_1}
\end{equation}
Since $\lambda$ consists always negative real part, so motion is asymptotically
stable in $z$ direction. Now  Eq.(\ref{eq:lambda_x}) and ({~\ref{eq:lambda_y}) has singular solution if,
\[
\begin{vmatrix}
\lambda^2-\lambda{U^0_{x\dot{x}}}-U^0_{xx}& -(2n+{U^0_{x\dot{y}}})\lambda-U^0_{xy} \\(2n-{U^0_{y\dot{x}}})\lambda-U^0_{yx}& \lambda^2-\lambda{U^0_{y\dot{y}}}-U^0_{yy}\\
\end{vmatrix}
=0
\]
\begin{equation}
\Rightarrow \quad \lambda^4+a\lambda^3+b\lambda^2+c\lambda+d=0 \label{eq:cheq}
\end{equation}
At the equilibrium points the system of Eq.(1),(2) gives
\begin{align}
a&=3\zws,\quad b=2n^2-f_*-\frac{3\mu{A_2}}{{r^5_2}_*}+\frac{2W^2_1}{{r^4_1}_*}=b_0+b_1 \quad  b_0=2n^2,b_1=-f_*-\frac{3\mu{A_2}}{{r^5_2}_*}+\frac{2W^2_1}{{r^4_1}_*}\\
c&=-a\biggl\{n^2+\frac{\mu{A_2}}{{r^5_2}_*}+\zd{\frac{\mu}{{{r^2_1}_*}{{r^5_2}_*}}\zabs{y^2_*}}\biggr\}\\
d&=(n^2-f_*)\biggl[n^2+2f_*-\frac{3\mu{A_2}}{{r^5_2}_*}\biggr]+\zd{\frac{9\mu\amc{q_1}}{{{r^5_1}_*}{{r^5_2}_*}}\zabs{y^2_*}}\\\notag
&-\zd{\frac{6nW_1}{{{r^4_1}_*}{{r^5_2}_*}}\zabs\Bigl\{\zxs^2-\zxs+y^2_*\Bigr\}}+\frac{n^2W^2_1}{{r^4_1}_*}
\end{align}
where \[{r_1}_*=q^{1/3}_1\left\{1-\frac{nW_1}{6\amc y_0}-\frac{A_2}{2} \right\},\quad {r_2}_*=1+\frac{nW_1}{3\mu y_0}\left(1-\frac{5}{2}A_2\right)\]
Let $q=1-\epsilon$ with $\epsilon>0$  very small real quantity then  \[f_*= 1-\frac{nW_1(9+2\epsilon)(4-3A_2)}{18\sqrt{3}}+3\frac{\amc A_2}{2}\]
We can write the four roots of the classical characteristic equation as,
\begin{equation}
\lambda_j=\pm zi,\quad \text{where},\quad z^2=\frac{1}{2}\Bigl\{1\mp\bigl[1-27\mu\amc\bigr]^{1/2}\Bigr\},\quad j=1,2,3,4.
\end{equation}
Due to P-R drag and oblateness,we assume the solution of the characteristic Eq.(~\ref{eq:cheq}) of the form
\[
\lambda=\lambda_j(1+e_1+e_2i)=\pm[-e_2+(1+e_1)i]z
\]
where $e_1,e_2$ are small real quantities. We have in first order approximation
\begin{gather*}
\lambda^2=[-(1+2e_1)-2e_2i]{z^2},\quad \lambda^3=\pm[3e_2-(1+3e_1)i]{z^3},\quad \lambda^4=[(1+4e_1)+4e_2i]z^4
\end{gather*}
Substituting these values in characteristic Eq.(~\ref{eq:cheq}) neglecting products of $e_1$ or $e_2$ with $a,b_1,c$, we get
\begin{align}
e_1&=\zd{\frac{-d+b_0{z^2}-z^4}{2z^2(2z^2-b_0)}}\\
e_2&=\zd{\frac{\mp{cz}\pm{az^3}}{2z^2(2z^2-b_0)}}
\end{align}
If $e_2\neq 0$,then the resulting motion of particle displaced from equilibrium points is asymptotically stable only when all the real parts of $\lambda$ are negative.
For stability, we require $Re(\lambda)<0$. Rate of growth or decay is determined by the largest real part of the roots of the characteristic equation, where  $Re(\lambda)$ denotes the real parts of $\lambda$ and given by

\begin{equation}
Re(\lambda)=\zd{\frac{c-az^2}{2(2z^2-1)}}
\end{equation}
we have,
\begin{equation}
z^2=\frac{1}{2}\Bigl\{1\mp\bigl[1-27\mu\amc\bigr]^{1/2}\Bigr\}
\end{equation}
Taking positive sign , \[ z^2=1-\frac{27}{4}\mu\amc  \]
Consider
\begin{equation*}
Re(\lambda)<0 \quad \text{then} 
\quad \Bigl\{c-a\bigl[1-\frac{27}{4}\mu\amc\bigr]\Bigr\}\Bigl\{1-\frac{27}{2}\mu\amc\Bigr\}^{-1}<0 
\end{equation*}
\begin{equation}
\Rightarrow \quad c+\frac{27}{4}\mu\amc(2c-a)<a
\end{equation}
Taking negative sign,
\begin{align}
Re(\lambda)&=-\Bigl[c-a\frac{27}{4}\mu\amc\Bigr]{\Bigl[1-\frac{27}{2}\mu\amc\Bigr]}^{-1}<0 \notag\\
&\Rightarrow  0<c+\frac{27}{4}\mu\amc(2c-a)\\
\text{From  above inequalities we get},\notag\\
&0<c+\frac{27}{4}\mu\amc(2c-a)<a\\
&0<c<a\quad \text{as}\quad \mu\rightarrow 0 \qquad \text{Murray C.D.(1994)}\label{eq:nc}
\end{align}
Inequality(~\ref{eq:nc}) is necessary condition for stability of triangular equilibrium points $L_4$ and $L_5$.

\section{Conclusion}

\begin{enumerate}
\item[(i)]  When $W_1=0$, then
\begin{align*}
a&=0,\- c=0,\- b=n^2-3\mu{A_2}, \\
d&=9\mu\amc\abc\biggl(1-\zd{\frac{q^{2/3}_1}{4n^4}}\biggr)n^2 
\end{align*}
Equation (~\ref{eq:cheq}) gives us,
\begin{align}
 \lambda^2&=\zd{\frac{-b\pm{(b^2-4d)}^{1/2}}{2}} \quad \text{For stable motion} \quad b^2>4d\notag  \\
\text{i.e.}, \quad
&{\bigl(n^2-3\mu{A_2}\bigr)^2}>36\zd{\mu\amc\abc\Bigl(1-\zd{\frac{q^{2/3}_1}{4n^4}}\Bigr)n^2}\notag
\end{align}
\item[(ii)] For classical case $A_2=0$, $q_1=1$, $W_1=0$, $n=1$. Then from inequality (~\ref{eq:nc}) we have
\[
1>27\mu\amc \quad \Rightarrow \quad \mu <.0385
\]
\item[(iii)] When $W_1\neq 0$, $A_2=0$, $n=1$.  Then
\begin{align*}
a&=3\zws, b=2-f_*+\frac{2W^2_1}{r^4_{1_*}},c=-a\Bigl\{1+\frac{\mu{y^2_*}}{{r^2_{1_*}}{r^5_{1_*}}}\Bigr\}\\
d&=(1-f_*)(1+2f_*)+\zd{\frac{9\mu\amc{q_1}}{r^2_{1_*}r^5_{1_*}}y^2_*}
-\zd{\frac{6W_1}{r^4_{1_*}r^5_{2_*}}\bigl[\zxs^2-\zxs+y^2_*\bigr]}+\frac{W^2_1}{r^4_{1_*}}
\end{align*}
as in  Chernikove Yu.A.(1970)\\
\end{enumerate}
 If we consider the P-R drag force then $c< 0$. This dose not satisfy the necessary condition of stability. Hence motion is unstable in the linear sense.
%\vspace{.1in}

{\bf Acknowledgment:}
We are thankful to D.S.T., Government of India, New Delhi for sanctioning a
project DST/MS/140/2K dated 02/01/2004 on this topic.\\

\end{document}